\documentclass[12pt]{amsart}
\usepackage{amssymb,amscd, fullpage}

\usepackage{palatino}
\usepackage{euler}

\numberwithin{equation}{section}

\newtheorem{theorema}{Theorem}
\newtheorem{theorem}{Theorem}[section]

\newtheorem{lemma}[theorem]{Lemma}
\theoremstyle{definition}
\newtheorem{definition}[theorem]{Definition}

\theoremstyle{remark}
\newtheorem{remark}[theorem]{Remark}

\newcommand{\R}{\mathbb{R}}
\newcommand{\C}{\mathbb{C}}
\newcommand{\Z}{\mathbb{Z}}

\newcommand{\T}{\mathbb{T}}
\newcommand\lie[1]{\mathfrak{#1}}

\newcommand{\fg}{\lie{g}}

\newcommand{\fk}{\lie{k}}

\def    \inv    {^{-1}}

\newcommand\cone{\mathaccent23c}

\newcommand{\bM}    {\overline{M}}
\newcommand{\intF} {{\mathaccent23{F}}}

\newcommand{\labell}[1]	{\label{#1}}

\begin{document}

\title{Toric  symplectic singular  spaces I: isolated singularities }

\author{D. Burns \and V. Guillemin \and E. Lerman}\thanks{Supported in
part by NSF grants DMS-0104047 (DB), DMS-0104116 (VG) and
DMS-0204448(EL)} 
\address{M.I.T., Cambridge, MA 02139 \and University of Michigan, Ann
Arbor, MI 48109}

\email{dburns@umich.edu}

\address{M.I.T., Cambridge, MA 02139}

\email{vwg@math.mit.edu}

\address{University of Illinois, Urbana, IL 61801 and \newline
\,\, Australian National University, Canberra, ACT 0200, Australia}

\email{lerman@math.uiuc.edu, lerman@maths.anu.edu.au}


\begin{abstract}
We generalize a theorem of Delzant classifying compact connected
symplectic manifolds with completely integrable torus actions to certain
singular symplectic spaces.  The assumption on
singularities is that if they are not finite quotient then they are
isolated.
\end{abstract}
\maketitle

\section{Introduction}

It is an old question as to what the symplectic geometry analogue of a
singular algebraic variety should be. For example, it is an exercise
in Gromov's book \cite{Gro}.  One would like to have a definition of a
singular object in symplectic geometry (not to be confused with
symplectic varieties in algebraic geometry) that includes, at the very
least, symplectic quotients by proper Lie group actions and, perhaps, affine
subvarieties of complex vector spaces $\C^n$.  There are a number of
reasons for wanting such a definition.
\begin{enumerate}
\item Symplectic analogues of algebro-geometric objects shed new light on 
questions of symplectic and algebraic geometry, providing new
understanding and new techniques.  For instance Alexeev and Brion
\cite{ABr} recently proved the existence of degenerations of spherical
varieties to toric varieties. The proof is not geometric; it uses
commutative algebra.  It would be interesting to find a purely
symplectic proof.  And for that one needs to understand what
symplectic singular spaces are.
 
\item  An intrinsic definition of singular symplectic 
spaces will play a useful role in understanding the dynamics of
symmetric Hamiltonian systems, in particular in understanding the
stability and bifurcation of relative equilibria and the equivariant
analogue of the Liapunov - Weinstein - Moser center theorem (see for
example \cite{GL} and references therein).
\item 
In the same spirit, a good definition of a singular symplectic space
should allow one to develop an analogue of Floer homology and to
start proving the existence of relative periodic orbits in symmetric
Hamiltonian systems without resort to perturbative methods, that is,
away from relative equilibria.
\item Such a definition should also shed new light on quantization of 
singular systems.
\end{enumerate}

\subsection*{What should a  symplectic singular  space be?}

Fourteen years ago Sjamaar and Lerman proved  \cite{SL} that
symplectic quotients by actions of compact Lie groups are stratified
spaces (the result is also true for proper actions of non-compact
groups, see \cite{BL} and \cite{GL}).  They also proved that
symplectic quotients have a number of interesting properties, some of
which were abstracted in a definition of a symplectic stratified space
\cite[Definition 1.12]{SL}:

\begin{definition} \label{def*}
A {\em symplectic stratified space} is a stratified space $X$ together with
a subalgebra $C^\infty (X)$ (``the algebra of smooth functions'') of
the algebra of continuous functions $C^0(X)$ such that
\begin{enumerate}
\item each stratum $S$ of $X$ is a symplectic manifold,
\item $C^\infty(X)$ is a Poisson algebra,  and
\item the embedding $S\hookrightarrow X$ is Poisson for each stratum $S$.
\end{enumerate}
\end{definition}
The last line means that for any function $f\in C^\infty(X)$ its
restriction $f|_S$ to any stratum $S$ is a smooth function in the
ordinary sense, and for any two functions $f, g\in C^\infty (X)$ we have
$
\{f|_S, g|_S\}_S = \{f, g\}|_S,
$ where on the left $\{\cdot, \cdot\}_S$ denotes the Poisson bracket
defined by the symplectic form on the stratum $S$.  This tells us that
a singular symplectic space should be a {\em stratified space}.
Recall that a stratified space is a topological space with a locally
finite decomposition into manifolds, which are called strata.
Moreover, stratified spaces are defined recursively.  Namely, a
neighborhood of a point is homeomorphic to a product of a ball with a
cone on a compact stratified space, called the link of the singularity
at the point.  Points in the same stratum have isomorphic links.

The results of \cite{SL} also seem to suggest that a singular
symplectic space should be some sort of a {\em differential space} in the
sense of Aronszajn
\cite{Arons}, Spallek \cite{Spallek} and Sikorski \cite{Sikorski} 
(see \cite{Navarro} for the state of the art on differential spaces
and a large bibliography).  We feel, however, that this is may not be
the best way to approach symplectic singular spaces.  For instance,
consider a complex algebraic subvariety $X$ of $\C^n$.  The most
natural space of smooth functions on $X$ is the space of restrictions
$C^\infty (\C^n) |_X$, the space of Whitney smooth functions.  But
there is no obvious Poisson bracket in sight.  Here is another
example.  Take $X =\R^2$ with the standard symplectic form and
stratify it as $X = (\R^2 \smallsetminus \{0\}) \cup \{0\}$.  Then for
any $n$, the space of $\Z/ n$-invariant functions on $\R^2$ is a
Poisson algebra, and embeddings of strata are Poisson.  In this
example there are infinitely many choices of $C^\infty (X)$ and no
preferred one.

In this paper we set ourselves a more modest goal.  We investigate
what a {\em toric} symplectic singular space (with isolated
singularities) should be.  The paradigm here is a theorem of Delzant
\cite{D}: compact symplectic manifolds with completely integrable
torus actions are smooth projective toric varieties.  The lesson here
is that toric symmetries turn a floppy symplectic manifold into a
rigid algebraic variety.

In the case of toric symplectic singular spaces the links of
singularities are toric {\em contact} spaces.  As the classification
of contact toric manifolds shows \cite{Lctm}, in order for the moment
map image to be locally convex, these contact toric spaces have to be
of Reeb type.  Recall that the notion of Reeb type was introduced by
Boyer and Galicki in their study of Sasakian toric manifolds
\cite{BoGa}.

Since symplectic orbifolds can be treated on the same footing as
symplectic manifolds \cite{LT} and since the category of orbifolds is stable
under symplectic cuts, we will assume that our spaces are stratified
by orbifolds, not manifolds.

The main result of the paper is Theorem~\ref{thm1}. We prove
that compact symplectic toric spaces, such that their non-orbifold
singularities are isolated, are classified by convex rational polytopes
that are simple away from the vertices.  Consequently they are
isomorphic to symplectic quotients of $\C^N$ by actions of compact
abelian Lie groups (i.e., products of tori $\R^n/\Z^n$ and of finite
abelian groups).  It then follows that all such spaces are projective
toric varieties.

\section{Classification of toric symplectic spaces with isolated 
 singularities}

In this section we define symplectic toric stratified spaces with the
property that their non-orbifold singularities are
isolated\footnote{Their orbifold singularities are arbitrary.} and
prove that they are classified by rational polytopes that are simple
away from vertices.  We start with a few preliminary definitions and
remarks.

One of the techniques of the proof involves blowing up singularities.
Now in the symplectic category blow-ups are symplectic cuts \cite{Lcuts}.
Symplectic cuts are symplectic quotients. So orbifold singularities
are generic and thus unavoidable. On the other hand, the symplectic
cuts of orbifolds are again, generically, orbifolds.  So it is more
convenient to work with orbifolds throughout.  Therefore we will
consider stratified spaces stratified by {\em
orbifolds}.\footnote{There is a fair amount of confusion in the
literature regarding the notions of orbifold and orbifold morphisms.
In this paper we will only deal with reduced (effective) orbifolds,
that is, orbifolds for which the structure groups all act effectively,
as defined, for example, by Moerdijk and Pronk \cite{MP}.  It seems to
us that the correct notion of morphism of effective orbifolds is that
of a {\em strong map} of Moerdijk and Prong ({\em op.\ cit.}).  In
many instances in symplectic geometry the notion of a reduced orbifold
is too restrictive.  The main reason is that slice representations for
locally free proper actions of Lie groups need not be effective.  The
general correct notion of orbifolds and their morphisms (``good
maps'') is probably the one due to Chen and Ruan \cite{CR}.  It is a
folklore wisdom that Chen and Ruan's definitions make orbifolds into
stacks over the category of orbifolds.  }  Thus, for us, a singular
space with isolated singularities is a Hausdorff topological space $M$
together with a discrete set of points $\{x_\alpha\} \subset M$ so
that $$ M_{reg}:= M\smallsetminus
\{x_\alpha\}
$$ 
is a smooth (meaning $C^\infty$) orbifold.  We also want the
neighborhoods of singularities $x_\alpha$ to have bounded 
topology.  Therefore we require that for each index $\alpha$ there is
a compact orbifold $L_\alpha$ (the link of the singularity at
$x_\alpha$ ) so that a neighorhood of $x_\alpha $ in $M$ is
homemorphic to the open cone $$
\cone (L_\alpha) : = \left([0, \infty)\times L_\alpha \right)/
 \left(\{0\} \times L_\alpha\right).
$$
In fact it will be more convenient for us to define the topological
cone $\cone(L)$ on a space $L$ by
$$
\cone (L) : = \left([-\infty, \infty)\times L \right)/
 \left(\{-\infty \} \times L\right), $$ and to refer to the image $*$
 of $\left(\{-\infty \} \times L\right)$ in $\cone (L)$ as the {\em
 vertex} of the cone.  Here we think of $[-\infty, \infty)$ as the
 topological space which is homeomorphic to $[0, \infty)$ and contains
 $\R = (-\infty, \infty)$ as a dense open subset.

By a torus we mean a compact connected abelian Lie group.  We denote
the Lie algebra of a torus $G$ by $\fg$.  We let $\fg^* = Hom (\fg,
\R)$, the dual of the Lie algebra, and $\Z_G = \ker\{\exp: \fg \to
G\}$, the integral lattice of $G$.  A vector $v\in \Z_G$ is {\em
primitive} if it is not a positive integer multiple of another vector
in the lattice $\Z_G$.  A codimension one face of a polytope is a {\em
facet}. A polytope in $\fg^*$ is rational if the normals of all of its
facets lie in the integral lattice $\Z_G$.  In general the normals to
the facets of a polytope in $\fg^*$ lie in $(\fg^*)^* =\fg$. Our
convention is that all the normals to the facets are primitive and
outward pointing.  A polytope in $\fg^*$ is {\em simple} if all the
facets are in general position.  In particular there are exactly $\dim
\fg^*$ edges coming out of every vertex.  Thus a cube is simple while
an octahedron is not. However, an octahedron {\em is} simple away
from the vertices.

A {\em symplectic toric $G$-orbifold} is a triple $(M, \omega, \Phi: M\to
\fg^*)$ where $M$ is a connected orbifold, $\omega$ is a symplectic form on $M$ and
$\Phi$ is a moment map for an effective completely integrable action
of a torus $G$ on $(M, \omega)$. We {\em don't} assume that our toric
orbifolds are compact, nor do we assume that the moment maps are
proper.

\begin{remark}\labell{rm-labels}
A classification theorem of Lerman and Tolman \cite{LT} asserts that
compact connected effective symplectic toric $G$-orbifolds are
classified by simple rational polytopes in $\fg^*$ together with
positive integers labels attached to each facet.  The need for the
integer labels can be seen in the following example.  Consider
the teardrop orbifold $X_n$, that is, consider an orbifold
whose underlying topological space is homeomorphic to the two sphere
and whose only singularity is modelled on $\C/ (\Z/n)$.  The orbifold
carries a symplectic form which is invariant under a Hamiltonian
action of the circle $S^1$.  The associated polytope is a closed
interval in $\fg^*\simeq \R$.  The classification theorem in question
tells us that the image of the singularity is one of the end points of
the interval and that this facet of the interval should be labelled
$n$.  More generally it turns out that for symplectic toric orbifolds,
the orbifold structure groups of points that the moment map sends to
facets are cyclic.  The integer labels are simply the orders of
these groups.
\end{remark}

For us a {\em symplectic cone} is a symplectic orbifold $(N, \omega)$
together with a free proper action $\rho_t$ of the reals such that
$\rho_t^* \omega = e^t\omega$ for all $t\in \R$.  Just as for
manifolds there is a one-to-one correspondence between symplectic
cones and (cooriented) contact orbifolds. If $(N, \omega, \rho_t)$ is
a symplectic cone, then $L = N/\R$ is a contact orbifold.  Conversely,
if $L$ is a contact cooriented orbifold then its symplectization is a
symplectic cone. These facts are standard for manifolds.  Their proofs
in the orbifold case carries over {\em mutatis mutandis} from the
manifold case: replace every occurance of the word ``manifold'' by
``orbifold.''

We will say that a symplectic cone $(N, \omega, \rho_t)$ is $G$-{\em
toric} if there is an effective action of a torus $G$ on $N$
preserving $\omega$ and commuting with dilations $\rho_t$ and having the
property that $2 \dim G = \dim N$.  Now, whenever a group action on a
symplectic cone preserves a symplectic form and commutes with
dilations $\rho_t$, it is Hamiltonian since it preserves the
contraction of $\omega$ with the vector field $\zeta$ generating the
dilations $\rho _t$.  Note that the associated moment map $\Phi$
defined by the primitive $\iota (\zeta) \omega $ is homogeneous: $\Phi
(\rho_t (x)) = e^t \Phi (x)$. Hence $\lim _{t \to -\infty} \Phi
(\rho_t (x)) = 0$, and $\Phi$ extends to a continous map on the
topological cone $\cone (N/\R) = \{*\} \cup N$.  Consequently {\em
any} moment map on a toric symplectic cone $N$ extends to a continuous
map on the cone $\cone (N/\R)$

A symplectic toric cone $(N, \omega, \rho_t, \Phi: N\to \fg^*)$ is of
{\em Reeb type } if there is a vector $X\in \fg$ so that $\langle \Phi
(x), X\rangle > 0$ for all $x\in N$. Here again $\Phi $ denotes a {\em
homogeneous} moment map.
 We will see below
(Lemma~\ref{lem??}) that if additionally the quotient $N/\R$ is compact
and connected then the moment map image $\Phi (N) \cup \{0\} \subset
\fg^*$ is a cone on a simple rational polytope.  If $N$ is a manifold,
the converse is true as well.  This follows from the classification of
contact toric manifolds \cite{Lctm}.  Note that in general, even if
$N/\R$ is compact, the cone $\Phi (N) \cup \{0\}$ need not be a {\em
convex} cone.   This happens, for example, when $N$ is the
symplectization of an overtwisted 3-sphere \cite{Lccuts}.

We are now in a position to define the main object of this section:
toric symplectic  singular 
spaces with isolated non-orbifold singularities.

\begin{definition} 
A {\em  toric symplectic singular  $G$-space with isolated
singularities} is a Hausdorff topological space $M$ with an effective
action of a torus $G$, which has a number of properties.  Namely, we
assume that
\begin{enumerate}
\item  There is a discrete set of points $\{x_\alpha\} \subset M$ 
such that $M_{reg}:= M \smallsetminus \{x_\alpha\} $ is a
symplectic toric $G$-orbifold.

\item Topologically the neighborhoods of the singular points $x_\alpha$ are 
cones on compact orbifolds, which are not quotients of odd dimensional 
spheres $S^{2n-1}$ by finite subgroups of $\T^n$.

\item Symplectically and equivariantly the punctured neighborhoods of singular 
points $x_\alpha $ are neighborhoods of $-\infty$ in toric symplectic
cones.  Hence by the preceding remarks the moment map $\Phi$ on
$M_{reg}$ extends to a continuous function on $M$ which we will
continue to denote by $\Phi$.

\item For each singular point $x_\alpha$ there is a vector 
$Y_\alpha \in \fg$ so that 
$\langle \Phi - \Phi (x_\alpha), Y_\alpha\rangle >0$ on $U_\alpha
\smallsetminus x_\alpha$, where $U_\alpha $ is a $G$-invariant
neighborhood of $x_\alpha$.  That is, the neighborhoods of
singularities are  symplectic cones of Reeb type.
\end{enumerate}
\end{definition}
The condition (2) on the links insures that the singularties are
genuine non-orbifold singularities.

\begin{theorema}\labell{thm1}
Let $(M, \omega, \Phi: M \to \fg^*)$ be a compact connected toric
symplectic singular $G$-space with isolated singularities.  Then
\begin{enumerate}
\item The moment image $\Delta: = \Phi (M) \subset \fg^*$ is a rational 
polytope, which is simple except possibly at vertices.

\item The fibers of the moment map $\Phi: M \to \fg^*$ are $G$-orbits.

\item The moment polytope $\Delta$ plus positive integers attached to
  the facets of $\Delta$ uniquely determine the topological space $M$, the
  singularities $\{x_\alpha\}$ and the symplectic toric $G$-orbifold
  $M_{reg} = M \smallsetminus \{x_\alpha\}$.
\end{enumerate}
\end{theorema}

\begin{remark}
\labell{rmk}
There is an easy converse to Theorem~\ref{thm1}.  Namely, given any
rational polytope $\Delta \subset \fg^*$ which is simple away from
vertices and a set of positive integers attached to facets, there
exists a singular symplectic $G$-variety with moment image $\Delta$.
The proof is essentially the same as in Delzant's original paper
\cite{D} modified as in \cite{LT} to take care of the positive
integer labels attached to facets.  We will say more about it later in
the paper.
\end{remark}
We now proceed with the proof of Theorem~\ref{thm1}.  We start by
recalling what symplectic cuts are \cite{Lcuts}.  Let $(M, \omega)$ be a
 symplectic
orbifold with an effective Hamiltonian action of a circle $T$ and an
associated moment map $\mu:M \to \R$.  Let $a$ be a regular value of
$\mu$.  Then the action of $T$ on the level set $\mu\inv (a)$ is
locally free and
$$
M/\!/_a T:= \mu \inv (a)/T
$$
is again naturally a symplectic orbifold by the reduction theorem of
Marsden, Weinstein and Meyer.  Consider now an orbifold with boundary
$$
 \{\mu \geq a\} := \{ x\in M \mid \mu (x) \geq a\}.
$$
Define an equivalence relation $\sim$ on $\{\mu \geq a\}$ by
identifying points on the boundary that lie on the same orbit of $T$.
The following theorem is one of the main observations of \cite{Lcuts}. We will
 not describe its proof here.
\begin{theorem}\labell{thm-cuts}
In the above notation, the topological space
$$
\overline{M}_{\geq a} := \{\mu \geq a\}/\sim
$$ can be naturally given the structure of a symplectic orbifold so
that the inclusion $M/\!/_a T \hookrightarrow \overline{M}_{\geq a}$ is
a symplectic embedding and the difference $\overline{M}_{\geq a}
\smallsetminus M/\!/_a T$ is symplectomorphic to the open subset 
$\{\mu > a\}$ of $(M, \omega)$. 
Moreover the restriction $\mu|_{\{\mu > a\}}$ descends to a moment map 
$\bar{\mu} : 
\overline{M}_{\geq a} \to \R$ for the action of $S^1$ on the cut space 
$\overline{M}_{\geq a}$.
\end{theorem}

\begin{remark}
 {\em Mutatis mutandis} the statement
 holds also for 
$$
\overline{M}_{\leq a} := \{\mu \leq a\}/\sim .
$$
\end{remark}
\begin{remark}\label{rm-local}
Theorem~\ref{thm-cuts} is {\em local}: If $\Sigma \subset (M, \omega)$
is a hypersuface which separates $M$ into two manifolds with boundary
$M^+$ and $M^-$ and if the null foliation of $\omega|_\Sigma$ consists
of closed leaves, then we may collapse the leaves of the null
foliation in the boundaries of $M^+$ and $M^-$ and obtain symplectic
orbifolds without boundary.
\end{remark}
\begin{remark}
The results of Theorem~\ref{thm-cuts} hold {\em equivariantly}. In
particular, suppose a torus $G$ acts on $(M, \omega)$ with moment map
$\Phi: M\to \fg^*$ and $T\subset G$ is a closed subgroup generated by
$X\in \fg$ so that $\mu = \langle \Phi , X\rangle$. Then the cut space
$\overline{M}_{\geq a}$ inherits a Hamiltonian action of $G$.
Moreover the restriction $\Phi|_{\{\mu > a\}}$ descends to a
$G$-moment map $\bar{\Phi}$ on the cut space and
$$
\bar{\Phi}(\overline{M}_{\geq a}) = \Phi (\{\mu > a\}) = \Phi (M)
\cap \{\eta \in \fg^* \mid \langle \eta, X \rangle \geq a\}.
$$
\end{remark}

Theorem~\ref{thm-cuts} was motivated by symplectic blowups.  Namely
let 
\[
(M,
\omega) = (\C ^n, \sqrt{-1} \sum dz_j \wedge \bar{z}_j)
\] 
and let $\mu (z) = |z|^2$.  Then $\overline{M}_{\geq a} $ is the
symplectic blow-up of $\C^n$ at the origin and $\overline{M}_{\leq a}$
is $\C P^n$.

Next we sort out neighborhoods of singular points.
\begin{lemma} \labell{lem??} 
Let $(N, \omega, \rho_t)$ be a symplectic cone with an action of a
torus $K$ which preserves the symplectic form and commutes with
dilations $\rho_t$.  Let $\mu:N \to \fk^*$ denote the associated
homogeneous moment map.  Suppose there is a vector $X\in \fk$ such
that
\begin{equation} \labell{eq1.1} 
\langle \mu (x), X \rangle > 0 \textrm{ for all } x \in N,
\end{equation}
that is, suppose $(N, \omega, \rho_t, \mu: N\to \fk^*)$ is of Reeb type.
If the base $B= N/\R$ is compact and connected then $\mu$ is proper as
a map into $\fk^* \smallsetminus \{0\}$, the fibers of $\mu$ are
connected and the moment cone $\mu (N)\cup \{0\}$ is a cone on a
convex rational polytope.
\end{lemma}

\begin{proof}
Since the base $B$ of the cone is compact, (\ref{eq1.1}) is an open condition.
Therefore we may assume that $T=\{\exp tX \mid
t\in \R\}$ is a circle.  Then $f = \langle \mu, X\rangle $ is the
moment map for the action of $T$ on $(N, \omega)$.  Recall that $\mu$
is defined by
\[
  \langle \mu , Y\rangle = \omega (\zeta, Y_N)  
\]
for all $Y\in \fk$, where $\zeta$ is the is the vector field
generating $\rho_t$ and $ Y_N$ denotes the vector field on $N$ induced
by $Y$. Consequently $\mu (\rho_t (x)) = e^t \mu (x)$ for all $t\in
\R$, $x\in M$.   Note also that since $f = \omega (\zeta, X_N) >0$ by
 assumption, $X_N\not = 0$ 
anywhere on $N$; hence the action of $T$ on $N$ is locally free.

We now argue that $f\inv (1)$ is diffeomorphic to $B$, hence is
compact and connected.  Since
\[
f (\rho_t (x)) = e^t f(x)
\]
the $\rho_t$-orbit map $\pi: f\inv (1) \to B$ is injective.  Moreover,
it is onto: pick a section $s:B\to N$ of $\pi:N\to B$ (any principal
$\R$-bundle is trivial so such a global section does exist).  Let
\[
\tau (x) =-\log (f(s(x))). 
\]
Then
\[
f (\rho_{\tau (x)} (s(x)))  = e^{\tau (x)} f (s(x)) = 1 .
\]

Note that $f\inv (1) \hookrightarrow N$ is a $K$-invariant section for
$\pi: N \to B$.  Hence $N$ is {\em $K$-equivariantly} diffeomorphic to
$B\times \R$.  Moreover under this identification of $N$ with $B\times
\R$, $f$ is simply the map $f(m,t) = e^t$. Therefore $f$ is proper as a map 
into $(0,\infty)$.  Hence $\mu:N\to \{ \eta \in \fk^* \mid \langle
\eta, X\rangle >0\}$ is proper as well. By \cite[Theorem 4.3]{LMTW}, 
$\mu (N)$ is convex and the fibers of $\mu$ are connected.  But more
can be said.

The restriction $\mu|_{f\inv (1)}$ descends to a moment map
$\overline{\mu}$ for the action of $K$ on the quotient $f\inv (1)/T$,
which is a compact symplectic orbifold.  By the orbifold version of
the Atiyah-Guillemin-Sternberg convexity theorem \cite{LT}
$\overline{\mu }( f\inv (1)/T)$ is a rational polytope.  Hence $\mu
(f\inv (1) ) = \overline{\mu }(f\inv (1)/T)$ is a rational polytope as
well.  Therefore $\mu (N) \cup \{0\} = \R^{>0} \mu (f\inv (1)) \cup
\{0\}$ is a cone on a rational polytope.
\end{proof}

\begin{proof}[Proof of Theorem~\ref{thm1}]
Since $M$ is compact, the set of singularities $\{x_\alpha\}$ is
finite.  By assumption a neighborhood $U_i$ of each singular point
$x_i$ is the neighborooh of $-\infty$ in a toric symplectic cone of
Reeb type.  In particular there is a vector $Y_i$ in the integral
lattice $\Z _G$ of $G$ such that $\langle \Phi - \Phi (x_i),
Y_i\rangle > 0$ on $U_i$.  By Lemma~\ref{lem??} the restriction of the
function $f_i = \langle \Phi , Y_i\rangle$ to $U_i$ is proper as a map
onto its image.  By the proof of the lemma  the action of the circle $G_i :=
\{\exp tY_i \mid t\in \R\}$ is locally free on $U_i \smallsetminus \{x_i\}$.

Now for every sufficiently small $\varepsilon >0$ consider 
\[
M(\varepsilon) := M \smallsetminus \bigcup_i \{ m\in U_i \mid f_i (m) <
f_i (x_i) + \varepsilon \}; 
\]
it is a compact symplectic orbifold with boundary. The connected
components of the boundary are the contact toric orbifolds $\{m\in U_i
\mid f_i (m) = \varepsilon + f_i (x_i)\}$.  Moreover by Theorem~\ref{thm-cuts} and
subsequent remarks, if we divide each component by the corresponding
group $G_i$, the various components of the boundary disappear and the
result is a compact symplectic toric orbifold $\bM (\varepsilon)$.  Furthermore,
$\Phi|_{M(\varepsilon)}$ descends to a moment map $\bar{\Phi}$ on $\bM
(\varepsilon)$ for the induced action of $G$.  By construction $\Phi
(M(\varepsilon)) = \bar{\Phi}(\bM (\varepsilon))$ and the fibers of
$\Phi|_{M(\varepsilon)}$ are $G$-orbits if and only if the fibers of
$\bar{\Phi}$ are.  Since $\bM (\varepsilon)$ is a symplectic toric
orbifold, $\bar{\Phi}(\bM (\varepsilon))$ is a simple rational polytope
\cite{LT}.  Moreover, the facets of $\bar{\Phi}(\bM (\varepsilon))$
that come from cuts have labels 1 attached to them (see
Remark~\ref{rm-labels}).
Therefore $\Phi (M(\varepsilon))$ is
a simple rational polytope and the fibers of $\Phi|_{M(\varepsilon)}$
are $G$-orbits.  Note that
$$
\bigcup _{\{\varepsilon > 0\}} M(\varepsilon) = M_{reg}.  
$$

We are now in a position to argue that $\Phi (M)$ is convex and that the
fibers of $\Phi :M \to \fg^*$ are $G$-orbits.  Suppose $a, b \in M_{reg}$
are two points.    Pick $\varepsilon$ small enough so that
$a, b \in M(\varepsilon)$.  Since $\Phi (M(\varepsilon)) \subset \Phi
(M_{reg})$ is convex, $[\Phi(a) \Phi (b)] \subset \Phi (M(\varepsilon))$.
Therefore $\Phi (M_{reg})$ is convex.  Also, if $a, b\in M_{reg}$ and $\Phi (a) =
\Phi (b)$ then we can choose $\varepsilon$ so that $a, b\in
M(\varepsilon)$. Hence $a \in G\cdot b$.  Therefore the fibers of
$\Phi: M \to \fg^*$ are $G$-orbits. We will need this fact in the proof of
uniqueness below.\\

Since $M_{reg}$ is dense in $M$, $\Phi (M)$ is contained in the closure
$\overline{ \Phi (M_{reg})}$.  On the other hand, since $M$ is compact,
$\Phi: M \to \fg^*$ is proper.  Hence $\Phi (M) =
\overline{\Phi(M_{reg})}$.  In particular $\Phi (M)$ is convex.\\

We now argue that $\Phi (M)$ is a polytope, that the points $\Phi (x_i)$ are
among the vertices of this polytope and that $\Phi\inv (\Phi (x_i)) =
\{x_i\}$ for all $i$.  We cut off the conical neighborhoods  of singularities.  That is, 
for sufficiently small $\varepsilon$ we have
$$ 
M = M (\varepsilon)\cup \bigsqcup_i \{
m\in U_i \mid f_i (m) \leq f_i (x_i) +
\varepsilon \} 
$$
We know that both $\Phi (M(\varepsilon))$ and 
$$
\Delta_i (\varepsilon):= \Phi(\{ f_i  \leq  f_i (x_i) + \varepsilon \} )
$$ are convex polytopes which have disjoint interiors (the latter is
true because the fibers of $\Phi: M_{reg}\to \fg^*$ are $G$-orbits).
Also, for each $i$, $\Delta _i (\varepsilon)$ intersects $\Phi
(M(\varepsilon))$ in a facet that doesn't contain the vertex $\Phi
(x_i)$.  Hence $\Phi (M)$ is a polytope, $\Phi (x_i)$ are vertices
(not necessarily {\em all } the vertices) and, since $\Phi|_{\{ m\in
U_i \mid f_i (m) \leq f_i (x_i) + \varepsilon \} }$ has the property
that the fiber above $\Phi (x_i)$ is the point $x_i$, we have that the
preimages of vertices of $\Phi (M)$ are single points.  We leave it to
the reader to check that $\Phi (M)$ is simple except possibly at the
vertices.

\subsection*{Uniqueness}

Suppose $(M_i, \omega_i, \Phi_i :M_i \to \fg^*, \{x_1 ^{(i)}, \ldots,
x_N ^{(i)}\})$, $i=1,2$ are two toric symplectic singular $G$-spaces.
Suppose $\Phi_1 (M_1) = \Phi_2 (M_2)$ and the associated integer
labels on the facets agree.  Since by assumption the links of
singularities are not finite quotients of odd dimensional spheres, a
vertex of the polytope $\Phi_i (M_i)$ is simple if and only if it is
the image of a point in $(M_i)_{reg}$.  Therefore, after some renumbering, 
 $\Phi_1 (x_i ^{(1)})
= \Phi_2 (x_i ^{(2)})$.  

The symplectic toric orbifolds $((M_i)_{reg},
\omega_i, \Phi_i)$, $i=1,2$, are {\em locally isomorphic } in the
sense of \cite[p. 4222]{LT}. That is, every point $p \in 
\Phi_1((M_1)_{reg}) = \Phi_2((M_2)_{reg})$ has an open neighborhood $U \subset
\Delta _{reg}$ and an equivariant diffeomorphism $\phi_U: \Phi_1^{-1}(U) \to
\Phi_2^{-1}(U)$ with $\phi_U^*\omega_2 = \omega_1$ and $\phi^* \Phi_2
= \Phi_1$ on $\Phi_1^{-1}(U)$.  Hence to prove uniqueness it is enough
to argue that the group $H^2(\Phi_i ((M_i)_{reg},
\Z_G)$ is trivial \cite[Section 7]{LT}.  But $\Phi_i ((M_i)_{reg})$ is a polytope 
with some vertices possibly deleted, hence is contractible.  Thus the
uniqueness follows and we get an equivariant diffeomorphism $\phi:
(M_1)_{reg} \to (M_2)_{reg}$ with $\phi^*\omega_2 = \omega_1$ and
$\phi^* \Phi_2 = \Phi_1$.  It remains to argue that $\phi$ extends to
a continuous map $\phi: M_1
\to M_2$ so that $\phi^* \Phi_2 = \Phi_1$ still holds.

Fix $i$.  Take any sequence $\{y_n\}$ in $(M_1)_{reg}$ converging to $x_i ^{(1)}$.  
Then
$$
\lim _{n\to \infty}\Phi_2 (\phi (y_n))= \lim _{n\to \infty} (\Phi_1 (y_n)) 
= \Phi _1 (x_i ^{(1)}) = \Phi_2 (x_i^{(2)}).
$$
Since $M_2$ is compact, we may assume that $\{\phi (y_n)\}$
converges.  Since $\Phi_2 (\lim_{n\to \infty} \phi (y_n)) = \Phi_2
(x_i^{(2)})$, and since $\Phi_2\inv (\Phi _2 (x_i ^{(2)})) = \{ x_i
^{(2)} \}$, $\lim _{n\to \infty} \phi (y_n) = x_i ^{(2)}$.  Since the
sequence $\{y_n\}$ is arbitrary, $\phi$ does extend to a continuous map
$\phi:M_1 \to M_2$.
\end{proof}

Theorem~\ref{thm1} has an easy converse which is a slight generalization of 
\cite[Theorem~8.2]{LT}.  Namely,
\begin{theorem}\labell{thm-existence}
Let $G$ be a torus.  Let $\fg$ denote its Lie algebra, and let
$\Z_G \subset \fg$ denote its integral lattice.
Given a rational polytope $\Delta \subset \fg^*$, which is simple except
 possibly at the vertices,
and a positive integer $m_\intF$ attached to each open facet $\intF$ 
of $\Delta$, 
there exists a compact toric symplectic singular $G$-space 
$(M, \omega, \Phi: M\to \fg^*, \{x_\alpha\})$  such that
$\Phi (M) = \Delta$ and the orbifold structure group of a point
in $M$ which maps to a open facet $\intF$ is $\Z/ m_\intF \Z$. 

Moreover, $(M,\omega,\Phi: M\to \fg^*, \{x_\alpha\})$ is a symplectic
quotient of $\C^N$ by a closed abelian subgroup of $SU(N)$, where $N$ is
the number of facets of $\Delta$ minus $\dim G$.
\end{theorem}
The proof is a straightforward modification of the proof of
\cite[Theorem~8.2]{LT}.  We omit it.

\begin{remark}
Combining the above theorem with the uniqueness part of
Theorem~\ref{thm1} we see that {\em every} toric symplectic singular
space is a symplectic quotient of $\C^N$.  Therefore, by
\cite[Theorem~2.8]{Sj}, every toric symplectic singular space is a
complex analytic space.  Furthermore by \cite[Lemma~2.16]{Sj} it is a
K\"ahler space in the sense of Grauert.  The K\"ahler structures on
these  spaces is described elsewhere \cite{BGL}.
\end{remark}

\end{document}